\documentclass[11pt]{amsart}
\usepackage[dvips]{graphicx} 
\usepackage{amssymb,amscd,color,latexsym,epsfig,xypic,hyperref}
\usepackage[mathscr]{euscript}
\usepackage{tikz-cd}
\usepackage{tikz}
\usetikzlibrary{calc,trees,positioning,arrows,chains,shapes.geometric,%
    decorations.pathreplacing,decorations.pathmorphing,shapes,%
    matrix,shapes.symbols}

\tikzset{
>=stealth',
  punktchain/.style={
    rectangle,
    rounded corners,
    draw=black, thick,
    minimum height=3em,
    text centered,
    on chain},
  line/.style={draw, thick, <-},
  element/.style={
    tape,
    top color=white,
    bottom color=blue!50!black!60!,
    minimum width=8em,
    draw=blue!40!black!90, very thick,
    text width=10em,
    minimum height=3.5em,
    text centered,
    on chain},
  every join/.style={->, thick,shorten >=1pt},
  decoration={brace},
  tuborg/.style={decorate},
  tubnode/.style={midway, right=2pt},
}

\DeclareFontFamily{OT1}{rsfs}{}
\DeclareFontShape{OT1}{rsfs}{n}{it}{<-> rsfs10}{}
\DeclareMathAlphabet{\curly}{OT1}{rsfs}{n}{it}

\DeclareFontFamily{U}{mathb}{\hyphenchar\font45}
\DeclareFontShape{U}{mathb}{m}{n}{
      <5> <6> <7> <8> <9> <10> gen * mathb
      <10.95> mathb10 <12> <14.4> <17.28> <20.74> <24.88> mathb12
      }{}
\DeclareSymbolFont{mathb}{U}{mathb}{m}{n}

\makeatletter
\newcommand{\eqnum}{\refstepcounter{equation}\textup{\tagform@{\theequation}}}
\makeatother

\renewcommand\;{\hspace{.6pt}}
\newcommand\C{\mathbb C}
\newcommand\Q{\mathbb Q}

\newcommand\Z{\mathbb Z}

\newcommand\Ed{E_\bullet}

\newcommand\LL{\mathbb L}

\newcommand\PP{\mathbb P}

\newcommand\TT{\mathbb T}

\renewcommand\P{\mathbf P}
\newcommand\F{\mathbf F}
\newcommand\E{\mathbf E}
\newcommand\M{\mathbf M}

\newcommand\cO{\mathcal O}
\newcommand\cA{\mathcal A}

\newcommand\cC{\mathcal C}

\renewcommand\({\big(}
\renewcommand\){\big)}

\newcommand\wt{\widetilde}

\makeatletter
\newcommand{\so}{\ \ext@arrow 0359\Rightarrowfill@{}{\hspace{3mm}}\ }
\makeatother
\newcommand{\rt}[1]{\xrightarrow{\ #1\ }}
\newcommand\To{\longrightarrow}

\newcommand\into{\hookrightarrow}
\newcommand\INTO{\ \ar@{^(->}[r]<-.2ex>}
\newcommand{\Into}{\,\ensuremath{\lhook\joinrel\relbar\joinrel\rightarrow}\,}
\newcommand\onto{\to\hspace{-3mm}\to}

\newcommand\Mapsto{\ \longmapsto\ }

\renewcommand\_{^{}_}

\newfont{\bigtimesfont}{cmsy10 scaled \magstep5}
\newcommand{\bigtimes}{\mathop{\lower0.9ex\hbox{\bigtimesfont\symbol2}}}
\renewcommand\={\ =\ }

\newcommand\udot{^{\bullet}}

\DeclareMathSymbol{\lefttorightarrow}{3}{mathb}{"FC}
\DeclareMathSymbol{\righttoleftarrow}{3}{mathb}{"FD}

\newcommand\rk{\operatorname{rank}}
\newcommand\vir{\operatorname{vir}}

\newcommand\loc{\operatorname{loc}}

\newcommand\coker{\operatorname{coker}}
\newcommand\im{\operatorname{im}}
\newcommand\id{\operatorname{id}}

\newcommand\Hom{\operatorname{Hom}}
\renewcommand\hom{\curly H\!om}

\newcommand\Ext{\operatorname{Ext}}

\newcommand\Spec{\operatorname{Spec}}

\newcommand\Sym{\operatorname{Sym}}

\newcommand\tot{\operatorname{tot}}

\newcommand\beq[1]{\begin{equation}\label{#1}}
\newcommand\eeq{\end{equation}}
\newcommand\beqa{\begin{eqnarray*}}
\newcommand\eeqa{\end{eqnarray*}}

\newcommand\arXiv[1]{\href{http://arxiv.org/abs/#1}{arXiv:#1.}}
\newcommand\mathAG[1]{\href{http://arxiv.org/abs/math/#1}{math.AG/#1.}}
\newcommand\alggeom[1]{\href{http://arxiv.org/abs/alg-geom/#1}{alg-geom/#1.}}
\newcommand\hepth[1]{\href{http://arxiv.org/abs/hep-th/#1}{hep-th/#1.}}

\DeclareRobustCommand{\SkipTocEntry}[3]{}
\makeatletter
\newcommand\@dotsep{4.5}
\def\@tocline#1#2#3#4#5#6#7{\relax
  \ifnum #1>\c@tocdepth 
  \else
    \par \addpenalty\@secpenalty\addvspace{#2}%
    \begingroup \hyphenpenalty\@M
    \@ifempty{#4}{%
      \@tempdima\csname r@tocindent\number#1\endcsname\relax
    }{%
      \@tempdima#4\relax
    }%
    \parindent\z@ \leftskip#3\relax \advance\leftskip\@tempdima\relax
    \rightskip\@pnumwidth plus1em \parfillskip-\@pnumwidth
    #5\leavevmode #6\relax
    \leaders\hbox{$\m@th
      \mkern \@dotsep mu\hbox{.}\mkern \@dotsep mu$}\hfill
    \hbox to\@pnumwidth{\@tocpagenum{#7}}\par
    \nobreak
    \endgroup
  \fi}
\makeatother

\makeatletter \@addtoreset{equation}{section} \makeatother
\renewcommand{\theequation}{\thesection.\arabic{equation}}

\newtheorem*{thm*}{Theorem}

\theoremstyle{definition}
\newtheorem{rmk}[equation]{Remark}

\title{Quantum Lefschetz without curves}
\author{Jeongseok Oh and Richard P. Thomas}

\begin{document}
\maketitle
\begin{abstract} \noindent
Given one quasi-smooth derived space cut out of another by a section of a 2-term complex of bundles, we give two formulae for its virtual cycle.

They are modelled on the the $p$-fields construction of Chang-Li and the Quantum Lefschetz principle, and recover these when applied to moduli spaces of (stable or quasi-) maps. When the complex is a single bundle we recover results of Kim-Kresch-Pantev.
\end{abstract}


\section*{Introduction}
The original Quantum Lefschetz principle \cite{Ko, Giv, Kim} applied to curves in a quintic threefold $Q\subset\PP^4$ cut out by a quintic equation $s_Q\in H^0(\cO_{\PP^4}(5))$. Let $\iota\colon M_Q\into M_{\PP^4}$ denote the moduli spaces of genus $g$ stable maps to $Q$ and $\PP^4$ in some fixed degree. When $g=0$ then $M_{\PP^4}$ carries a natural bundle $E$ and section $s\in H^0(E)$ which, over the point $(f\colon C\to\PP^4)$, have fibre
\beq{sq}
f^*s_Q\ \in\ H^0\(f^*\cO_{\PP^4}(5)\).
\eeq
Clearly \eqref{sq} vanishes if and only if $\im f\subset Q$, so $s$ cuts $M_Q$ out of $M_{\PP^4}$ set-theoretically. In fact this is true as \emph{Deligne-Mumford stacks with perfect obstruction theory}: $s^{-1}(0)$ inherits a natural perfect obstruction theory which agrees with the one on $M_Q$. As a result its virtual cycle can be computed in terms of data on $M_{\PP^4}$:
\beq{QL}
\iota_*[M_Q]^{\vir}\=e(E)\cap[M_{\PP^4}].
\eeq
This aids computation because $\PP^4$ carries a torus action ($Q$ does not).

When $g\ge1$ the bundle $E$ is replaced by a \emph{2-term complex of vector bundles} $\Ed$ over $M_{\PP^4}$, which at the point $(f\colon C\to\PP^4)$ computes $H^*(f^*\cO_{\PP^4}(5))$. Motivated by work of Guffin-Sharpe-Witten \cite{GS}, Chang-Li \cite{CL1} moved the problematic $H^1$ term\footnote{Note this term in degree 2 is surjected onto by the obstruction sheaf of $M_{\PP^4}$, so the virtual tangent bundle of $M_Q$ is still supported in only degrees 0,\,1.} from degree 2 in the virtual tangent bundle of $(f^*s_Q)^{-1}(0)$ to degree 0 and dualised, forming a moduli space $F$ of stable maps to $\PP^4$ with ``$p$-fields". This is a cone over $M_{\PP^4}$ whose fibre over $(f\colon C\to\PP^4)$ is
\beq{pfield}
H^0\(\omega_C\otimes f^*\cO_{\PP^4}(-5)\)\ \cong\ H^1\(f^*\cO_{\PP^4}(5)\)^*.
\eeq
It comes with a natural perfect obstruction theory on which $f^*s_Q$ induces a natural \emph{cosection} \cite{KLi} whose zero locus is $M_Q\subset M_{\PP^4}\subset F$. Thus by cosection localisation \cite{KLi} we get a virtual cycle $[F]^{\loc}\in A_0(M_Q)$ localised to $M_Q$, which Chang-Li show recovers $M_Q$'s natural virtual cycle,
$$
[M_Q]^{\vir}\=(-1)^e\;[F]^{\loc}\ \in\ A_0(M_Q),
$$
where $e:=\rk\Ed$. Although they did not give an analogue of \eqref{QL} (but see the result \eqref{euler} below), we have gained something: by expressing $[M_Q]^{\vir}$ in terms of data on $M_{\PP^4}$ we can apply torus localisation. This has lead to great progress \cite{CGL, CL2, GJR, KLh} in computing higher genus Gromov-Witten invariants of the quintic.

Others \cite{BN,CG,CJW,CL3,Ke,KO,Lee,Pi1} have generalised this construction, but always for moduli of (stable or quasi-) maps of curves to a variety. We assumed this restriction was required to produce the cosection, but it turns out to exist more generally. 
The general setup replaces $M_Q\subset M_P,\,\Ed,\,s,\,F/M_P$ by data $\M\subset\P,\,\Ed,\,s,\,\F/\P$ as follows. We fix
\begin{itemize}
\item a quasi-smooth\footnote{A quasi-smooth derived structure $\bf P$ on an underlying scheme or stack $P$ is a slight enhancement of a perfect obstruction theory on $P$ which always exists in nature and seems to be necessary to set up our problem correctly.} ambient derived Deligne-Mumford stack $\P$ of (virtual) dimension $p$, whose underlying stack $P:=\pi_0({\bf P})$ has the resolution property (see \cite[Proposition 5.1]{Kr} for equivalent conditions),
\item an object $\Ed\in D(\mathrm{coh}\,\P)$ of rank $e$, quasi-isomorphic on $P$ to a 2-term complex of vector bundles $E_0\to E_1$,
\item a section $s\in H^0(\Ed)$ inducing a derived structure of dimension $p-e$ on its zero locus $\iota\colon\M:=s^{-1}(0)\into\P$.
\end{itemize}

Let $\F:=\Spec\Sym\udot(\Ed[1])\rt\pi\P$ be the quasi-smooth total space of the derived dual $\Ed^\vee[-1]$. The fibre of the underlying stack $F\subset\F$ over a point $x\in P$ is $H^1\((\Ed)_x\)^*$, just as in \eqref{pfield}. The section $s\in H^0(\Ed)$ \mbox{is equivalent} to a shifted function $\tilde s\in H^0\(\cO_{\F}[-1]\)$ on $\F$, linear on the fibres, via
$$
H^0(\Ed)\ \subset\ H^0\(\Sym\udot(\Ed[1])[-1]\)\=H^0\(\pi_*\;\cO_{\F}[-1]\)\=H^0\(\cO_{\F}[-1]\).
$$
Its derivative gives a map $d\;\tilde s\colon\TT_{\F}\to\cO_{\F}[-1]$. Taking $h^1$ therefore gives a map from the obstruction sheaf of $\F$ to its functions, i.e. a cosection \cite{KLi}.

\begin{thm*} The cosection $h^1(d\tilde s)\colon h^1(\TT_{\F})|_F\to\cO_F$ has image the ideal sheaf of $M\subseteq F$ if and only if $\M$ is quasi-smooth. In this case its cosection-localised virtual cycle is
\beq{thm}
[\F]^{\loc}\=(-1)^e\;[\M]^{\vir}\ \in\ A_{\;p-e}(M,\Q)
\eeq
Moreover we have the following analogue of \eqref{QL},
\beq{euler}
\iota_*\(e(E_1)\cap[\M]^{\vir}\)\=e(E_0)\cap[\P]^{\vir}\ \in\ A_{\;p-e_0}(P).
\eeq
\end{thm*}
\smallskip

When $\Ed=E_0$ is a bundle the formula \eqref{euler} is a result of Kim-Kresch-Pantev \cite{KKP}. In general  \eqref{euler} is the best we can do\,---\,we cannot expect a formula for the bare $[\M]^{\vir}$ in terms of $[\P]^{\vir}$ since the latter will have smaller dimension than the former when $\Ed=E_1[-1]$.

The first part of the Theorem describes $M$ as the critical locus of $\tilde s$, giving it a derived structure of dimension $2(p-e)$: \emph{twice} that of $\M$. In fact this derived structure on $M$ is the $(-2)$-shifted cotangent bundle $T^*[-2]\;\M$. In Section \ref{idea} we explain that the main idea behind the equality \eqref{thm} is that $T^*[-2]\;\M$ admits a virtual cycle \cite{OT1} computed using an auxiliary maximal isotropic subbundle of a certain bundle with quadratic form. Choosing one such subbundle gives $[\M]^{\vir}$, using another leads to $(-1)^e\;[\F]^{\loc}$.

It should perhaps not be a surprise that $(-2)$-shifted symplectic geometries \cite{PTVV}\,---\,and their associated virtual cycles \cite{BJ, OT1}\,---\,should play a role, given their relationship to cosections \cite{KP, Pi2} and to closed $(-1)$-shifted 1-forms and shifted critical loci. For this and other reasons this paper is in many ways just the $(-1)$-shift of the paper \cite{JT}.

We can apply our results to curve counting to recover the results of \cite{CJW, Pi1}. Let $M_P$ be the moduli space of stable maps of fixed degree and genus to a smooth projective variety $P$, with universal curve and map
$$
\xymatrix@R=15pt{\cC \ar[d]_\pi\ar[r]^f& P \\ M_P.\!\!}
$$
Let $(E,s)$ be a bundle and regular section over $P$ with smooth zero locus $Q\subset P$. Let $(\Ed\;,s)$ on $M_P$ be defined by the composition
$$
\cO_{M_P}\rt{\pi^*}R\pi_*\cO_{\cC}\rt{R\pi_*f^*s}R\pi_*f^*E\ =:\ \Ed\,.
$$
Then $s$ has zero locus $M_Q\subset M_P$ the stable maps to $Q$ and \eqref{thm} applies. There is a similar story for quasimaps but with $(E,s)$ defined on Artin quotient stacks $Q\subset P$. The same construction induces $(\Ed,s)$ cutting $M_Q$ out of the quasi-smooth Deligne-Mumford stack $M_P$, so again \eqref{thm} applies. 

\subsection*{Plan of paper} Section \ref{Euler} proves \eqref{euler}, while Section \ref{cosec} shows the cosection $h^1(d\tilde s)$ cuts out $M\subset F$. This leaves the proof of \eqref{thm} to Section \ref{idea} (if the local Kuranishi model can be globalised) and Section \ref{last} (in general).

By now both ``quantum" and ``Lefschetz" are both sufficiently far from our results that they should be thought of merely as motivation; ``virtual Euler" might be more appropriate.

We denote the derived dual $R\hom(E,\cO)$ of an object $E$ by  $E^\vee$, but use $E^*$ in the special case of vector spaces and bundles.

We use $\Q$ coefficients for our Chow groups throughout, but it should be noted that when $P$ is a scheme all results hold with $\Z$ coefficients. 
This uses the existence of the global maximal isotropic subbundle $\Lambda$ \eqref{cf} in Section \ref{last}, which ensures the results of \cite{OT1} hold with $\Z$ coefficients; see for instance \cite[Equation (34)]{OT1}. (In general \cite{OT1} works with $\Z[\frac 12]$-coefficients\,---\,results are proved on a certain bundle over $M$ on which a $\Lambda$ exists; inverting 2 then allows the descent back down to $M$.)

\medskip\noindent\textbf{Acknowledgements.} This paper is dedicated to Bumsig Kim, whose influence is all over this topic.

We thank Young-Hoon Kiem, Hyeonjun Park and Renata Picciotto for discussions about cosections, $p$-fields and their closely related preprints \cite{KP, Pi2}, which give a generalisation of our result to the non-quasi-smooth case.

Thanks also to Felix Janda and Bhamidi Sreedhar for useful conversations.

Both authors were supported by funding from a Royal Society research professorship.

\section{Euler classes}\label{Euler}
In this section we show \eqref{euler} follows easily from work of Kim-Kresch-Pantev. Let $s_0$ denote the projection of $s$ under $\Ed\to E_0$. This gives $M$ a different quasi-smooth structure\footnote{In local Kuranishi models for these structures, like in \eqref{Pfn} below, the cdgas for $\M$ and $\M'$ differ only in degrees $-2$ (by an $E_1^*$ term) and lower, so $\pi_0(\M')=\pi_0(\M)=M$.} $\M'$, cut out of $\P$ by $s_0\in H^0(E_0)$. By \cite{KKP} its virtual cycle pushes forwards to the right hand side of \eqref{euler},
$$
\iota_*[\M']^{\vir}\ =\ e(E_0)\cap[\P]^{\vir}.
$$
Let $\E:=\Spec\Sym\udot\Ed^\vee$ denote the total space of $\Ed$ with zero section $0_{\;\E}$ and section $s$, so that $\M$ is defined by the derived fibre product
$$
\xymatrix@R=15pt{
\M \ar[r]\ar[d]& \P \ar[d]^s \\
\P \ar[r]^{0_{\;\E}}& \E.\!\!}
$$
This gives the central horizontal exact triangle of the following commutative diagram of tangent complexes,
\beq{fp}
\xymatrix@R=15pt@C=35pt{
& \TT_{\P}\big|_{\M} \ar[d]_{(1,0)}\ar[r]^{0_{\;\E*}}& \TT_{\E}\big|_{\M} \ar@{=}[d] \\
\TT_{\M} \ar[r]\ar@{=}[d]& (\TT_{\P}\oplus\TT_{\P})\big|_{\M} \ar[r]^-{(0_{\;\E*},s_*)}\ar[d]_{(0,1)}& \TT_{\E}\big|_{\M} \\
\TT_{\M} \ar[r]& \TT_{\P}\big|_{\M}.}
\eeq
The zero section $0_{\E}\colon\P\into\E$ and the projection $\E\to\P$ together define a splitting of the tangent complex of $\E$ restricted to $\P$,
$$
\TT_{\E}|_{\P}\ \cong\ \TT_{\P}\,\oplus\,\TT_{\E/\P}|_{\P}\ \cong\ \TT_{\P}\,\oplus\,\Ed.
$$
So, up to shifts, the cone on the top row of \eqref{fp} is $\Ed|_{\M}$, while the bottom row is $\TT_{\M/\P}$. Since the central row is exact the upshot is that $\TT_{\M/\P}\cong\Ed|_{\M}[-1]$ sat in the exact triangle
$$
\TT_{\M}\To\TT_{\P}|_{\M}\To\TT_{\M/\P}[1]\,=\,\Ed.
$$
Restricting this to $M$ gives the top row of the following diagram of exact triangles; repeating the working with $(\Ed,s)$ replaced by $(E_0,s_0)$ gives the second row,
$$
\xymatrix@R=12pt{
\TT_{\M}\big|_M \ar[d]\ar[r]& \TT_{\P}\big|_M \ar@{=}[d]\ar[r]& \Ed\big|_M \ar[d] \\
\TT_{\M'}\big|_M \ar[d]\ar[r]& \TT_{\P}\big|_M \ar[r]& E_0\big|_M \ar[d] \\
E_1\big|_M[-1] && E_1\big|_M\,.\!\!\!}
$$
Picking a global locally free resolution $A\to K'$ for $T_{\M'}|_M$, the composition $K'\to\TT_{\M'}|_M[1]\to E_1|_M$ is onto, defining an exact sequence
$$
0\To K\To K'\To E_1|_M\To0
$$
and a resolution $A\to K$ for $\TT_{\M}|_M$.
Dualising, $\{K^*\to A^*\}=\LL_{\M}|_M\to\LL_M$ is a perfect obstruction theory for $M$ and the construction of Behrend-Fantechi \cite{BF} gives a cone $C\subset K$ such that $[\M]^{\vir}=0^{\;!}_{K}[C]$. Then
$$
[\M']^{\vir}\=0^{\;!}_{K'}[C]\=e(E_1|_M)\cap0^{\;!}_K[C]\=e(E_1)\cap[\M]^{\vir},
$$
the second equality by \cite[Theorem 6.3]{Fu}. This is the left hand side of \eqref{euler}.

\section{Cosection}\label{cosec}
Throughout it will be convenient to extend our $\,\wt{}\,$ notation: given any object $G\in D(\mathrm{coh}\,\P)$, a section $\varphi\in H^0(\P,\,G\otimes\Ed)$ can be thought of as a shifted section $\tilde\varphi\in H^0(\F,\,\pi^*G[-1])$, linear on the fibres of $\pi$, via
\begin{eqnarray}\nonumber
\varphi &\in& H^0(G\otimes\Ed)\ \subset\ H^0\(G\otimes\Sym\udot(\Ed[1])[-1]\) \\ \label{til}
&=&H^0\(G\otimes\pi_*\;\cO_{\F}[-1]\)\=H^0\(\pi^*G[-1]\)\ \ni\ \tilde\varphi.
\end{eqnarray}

In this Section we prove that $h^1(d\tilde s)$ cuts out $M\subseteq F$ if $\M$ is quasi-smooth. It is enough to work locally, where we can put everything in a standard model.

$\bullet$ The local model for a quasi-smooth $\P$ is a \emph{Kuranishi chart} $(A,B,t)$: a smooth ambient space $A$ over which we have a section $t$ of a bundle $B$ cutting out $\P$ in the sense that its structure sheaf is (quasi-isomorphic as a cdga to) the Koszul complex
\beq{Pfn}
\cO_{\P}\ \cong\ (\Lambda\udot B^*,t)\ :=\ \big\{\cdots\To\Lambda^2B^*\rt{t}B^*\rt{t}\cO_A\big\}.
\eeq
The cotangent complex of $\P$ is most easily described when $A$ is sufficiently small that $B$ admits a connection $D$. Then
\beq{LLP}
\LL_{\P}\=\(B^*\rt{\cdot Dt}\Omega_A\)\otimes\_{\cO_A}
(\Lambda\udot B^*,t)
\eeq
is a differential graded module over the dga $(\Lambda\udot B^*,t)$ in the obvious way. The exterior derivative acts on degree 0 functions by $\cO_A\ni f\mapsto df\otimes 1\in\Omega_A\otimes\cO_A$ and on degree $(-1)$ functions by
\beq{df}
B^*\,\in\,f\,\Mapsto\,(f\otimes1)\,\oplus\,Df\ \in\ (B^*\otimes\cO_A)\,\oplus\,(\Omega_A\otimes B^*).
\eeq
This can be checked to intertwine the differentials on the dga $(\Lambda\udot B^*,t)$ and the dgm \eqref{LLP}, and extends to degree $\le(-2)$ functions by the Leibniz rule.

$\bullet$ The local model for $\Ed$ is a complex of bundles $d\colon E_0\to E_1$ over $A$ tensored over $\cO_A$ with \eqref{Pfn}\,---\,i.e. the total complex $\Lambda\udot B^*\otimes\{E_0\to E_1\}$ with differential $t\otimes1-(-1)^i\otimes d$ on $\Lambda^i B^*\otimes E_0$. In degrees 0 and 1 this is
$$
\cdots\To(B^*\otimes E_1)\,\oplus\,E_0\xrightarrow{\,t-d\ }E_1,
$$
so $s\in H^0(\Ed)$ is represented by $(s_1,s_0)\in\Gamma\((B^*\otimes E_1)\oplus E_0\)$ such that $(t-d)(s_1,s_0)=0$. Thus $s_1(t)=d\circ s_0$ and $s$ is represented by a commutative diagram of $\cO_A$-modules
\beq{localmodel}
\xymatrix@R=18pt{
\cO_A \ar[r]^{s_0}\ar[d]_t& E_0 \ar[d]^d \\
B \ar[r]^{s_1}& E_1.\!\!}
\eeq

$\bullet$ The local model for $\F$ is inside $\tot\_{\!A}(E_1^*)$, cut out by the section
$$
r\ :=\ \(p^*t,\,-\wt{d^*}\)\ \in\ \Gamma\(p^*B\oplus p^*E^*_0\).
$$
Here $p$ is the projection $\tot\_{\!A}(E_1^*)\to A$ on which we are using the $\,\wt{}\,$ notation of \eqref{til}. Thus $\cO_{\F}$ is the associated Koszul cdga $\(\Lambda\udot(p^*B^*\oplus p^*E _0),\,r\)$ on $\tot\_{\!A}(E_1^*)$. In its degree $(-1)$ piece lies the $(-1)$-shifted function $\tilde s$,
$$
\tilde{s}\=\tilde s_1+p^*s_0\ \in\ \Gamma\(p^*B^*\oplus p^*E_0\).
$$
By \eqref{df} we can read off its exterior derivative $d\tilde s$ in
$$
\LL_{\F}\=\Big(p^*B^*\oplus p^*E_0\To\Omega_{\mathrm{tot}\_{\!A}(E_1^*)}\Big)\otimes\Lambda\udot(p^*B^*\oplus p^*E_0).
$$
It lies in the degree $(-1)$ part
$$
(p^*B^*\otimes\cO)\ \oplus\ (p^*E_0\otimes\cO)\ \oplus\ (\Omega\otimes p^*B^*)\ \oplus\ (\Omega\otimes p^*E_0),
$$
with respect to which it is
$$
\(\tilde s_1\otimes 1,\ 
p^*s_0\otimes1,\ \wt{Ds_1},\ p^*Ds_0\).
$$

Restricting to $F\subset\F$ kills the third and fourth terms (since they involve degree $(-1)$ functions). So we are left with showing that
\beq{cosex}
(\tilde s_1,p^*s_0)\,\colon\,p^*(B\oplus E^*_0)\To\cO_F
\eeq
has image the ideal sheaf of $M$. But $s_0$ cuts out $M$ from $P$ so $p^*(s_0)$ cuts out $F_M:=F\times_PM$ and what remains to show is that $\tilde s_1\colon p^*B|_{F_M}\to\cO_{F_M}$ generates the ideal sheaf of $M\subset F_M$. For this we consider the diagram
$$
\xymatrix@C=30pt{
p^*(B\oplus E_0)\big|_{F_M} \ar@/_1pc/[rr]_{p^*(s_1,-d)}\ar[r]^(.6){(1,0)}& p^*B\big|_{F_M} \ar[r]^-{p^*(s_1)}\ar`u[r]`[rr]^{\tilde s_1}[rr]+<0ex,1.6ex>& p^*E_1\big|_{F_M} \ar[r]_-{\tau|\_{F_M}}& \cO_{F_M},}
$$
where $\tau\in H^0(p^*E_1^*)$ is the tautological section of $E_1^*$ on $\tot\_{\!A}(E_1^*)$. Since $\tau\circ p^*(s_1)=\tilde s_1$ and $\tau\circ p^*(d\;)=\tilde d$\,---\,and the latter vanishes on $F_M$ because $\tilde d^*\colon p^*E_1^*\to p^*E_0^*$ vanishes on $F$\,---\,it follows that the diagram commutes. Most importantly, the image of $\tau\colon p^*E_1\to\cO_{\tot\_{\!A}(E_1^*)}$ is the ideal of the zero section $A\subset\tot\_{\!A}(E_1^*)$.
Now $\M$ is quasi-smooth if and only if
$$
h^2(\TT_{\M})\=0\=\coker\big[(s_1,-d\;)\colon(B\oplus E_0)\big|_M\to E_1\big|_M\big];
$$
see \eqref{TTM} below, for instance. So in this case the lower curved arrow is onto and its composition with $\tau|_{F_M}$ generates the ideal of the zero section $M\subset F|_M$. Thus so does the upper curved arrow, as required.


\section{Idea of the proof}\label{idea} We explain how \eqref{thm} works in the special case that \emph{the local model \eqref{localmodel} holds globally.} The main idea is to consider a \emph{third} derived structure on $M$, different from both $\M$ and $\M'$, namely the $(-2)$-shifted cotangent bundle $T^*[-2]\;\M$. This has a virtual cycle constructed in \cite{OT1} using a choice of maximal isotropic subbundle of a certain orthogonal bundle. Using one choice will recover $[\M]^{\vir}$, using another naturally gives $(-1)^e\;[\F]^{\loc}$. \medskip


From the model \eqref{localmodel}, $M\subset A$ is cut out by the section $(s_0,t)$ of the complex $E_0\oplus B\to E_1$ with differential $(-d,s_1)$. This endows it with a derived structure $\M$ with structure sheaf the Koszul complex
$$
\cO_{\M}\ \cong\ \Big(\Sym\udot\!\big\{E_1^*\xrightarrow{}E_0^*\oplus B^*\big\},(s_0,t)\Big).
$$
In particular its tangent complex $\TT_{\M}|_M$ is
\beq{TTM}
T_A\big|_M\xrightarrow{D(t,s_0)}B\big|_M\oplus E_0\big|_M\xrightarrow{(s_1,-d)}E_1\big|_M\,
\eeq
so $\M$ is quasi-smooth if and only if $(s_1,-d)\colon(B\oplus E_0)|_M\to E_1|_M$ is onto. In this case we set $K$ to be its kernel, so that
$$
\TT_{\M}\big|_M\=\big\{T_A\big|_M\xrightarrow{D(t,s_0)}K\big\}.
$$
Dualising induces a perfect obstruction theory $\LL_{\M}|_M\to\LL_M$, yielding a Behrend-Fantechi virtual cycle $[\M]^{\vir}$.\medskip

The shifted cotangent bundle $T^*[-2]\;\M$ has the same underlying stack $M$ but a different derived structure, with tangent complex
\beq{TT}
\TT_{T^*[-2]\;\M}\big|_M\=\big\{T_A\big|_M\xrightarrow{D(t,s_0)\,\oplus\,0\,}K\oplus K^*\xrightarrow{0\,\oplus\,D(t,s_0)^*}\Omega_A\big|_M\big\}.
\eeq
Since $T^*[-2]\;\M$ is $(-2)$-shifted symplectic it also admits a virtual cycle \cite{OT1}.
This depends on a choice of orientation (in the sense of \cite[Section 2]{OT1}) on the orthogonal bundle $K\oplus K^*$\,---\,making it an $SO(2k,\C)$ bundle\,---\,and the construction involves picking a maximal isotropic subbundle (though the final result is independent of it).

There is a canonical orientation $o_K$ on $K\oplus K^*$ which makes $K\subset K\oplus K^*$ a \emph{positive} maximal isotropic subbundle \cite[Equation (18)]{OT1}. Then picking $K^*\subset K\oplus K^*$ as our maximal isotropic subbundle, $T^*[-2]\;\M$'s virtual cycle is $[\M]^{\vir}$ by \cite[Section 8]{OT1}.\medskip

Applying the construction of \cite{OT1} to a different maximal isotropic, however, will lead naturally to the space $\F$.
We begin by replacing \eqref{TT} by the quasi-isomorphic complex
\beq{t*t}
(T_A\oplus E_1^*)\big|_M\To
(B\oplus E_0\oplus E_0^*\oplus B^*)\big|_M \To (\Omega_A\oplus E_1)\big|_M.
\eeq
Here the first arrow is the direct sum of $D(t,s_0)\colon T_A|_M\to(B\oplus E_0)|_M$ and $(-d^*,s_1^*)\colon E_1^*|_M\to(E_0^*\oplus B^*)|_M$, and the second arrow is its dual.
We claim \eqref{t*t} is the tangent bundle of the $(-2)$-shifted symplectic derived Deligne-Mumford stack cut out of $p\colon\tot\_{\!A}(E_1^*)\to A$ by the \emph{isotropic} section
\beq{sig}
\sigma\ :=\ \(p^*(s_0),\,p^*(t),\,-\tilde d^*,\,\tilde s_1^*\)\ \text{ of }\ p^*(B\oplus E_0\oplus E_0^*\oplus B^*).
\eeq
Here we use the $\,\wt{}\,$ notation of \eqref{til} and the natural quadratic form $q$ on $B\oplus E_0\oplus E_0^*\oplus B^*$\,---\,pairing $B\oplus E_0$ with its dual\,---\,so the commutativity of \eqref{localmodel} makes $\sigma$ isotropic. Thus we get a ``Darboux chart"
\beq{Darboux}
\xymatrix@=0pt{
& \(p^*(B\oplus E_0\oplus E_0^*\oplus B^*),\,p^*q\)\ddto \\  && \hspace{1cm}p^*q\;(\sigma,\sigma)=0,\hspace{-5mm} \\
 M\ =\ \sigma^{-1}(0)\ \subset\hspace{-17mm} & \mathrm{tot}\_{\!A}(E_1^*),\ar@/^{-2ex}/[uu]_(.45)\sigma}
\eeq
such that the two arrows of \eqref{t*t} are $D\sigma$ on $M$, which proves our claim. \medskip

The key observation is the following. Let $\Lambda$ denote the maximal isotropic subbundle $p^*(B^*\oplus E_0)$ and split \eqref{sig} as $\sigma=(\sigma_1,\sigma_2)\in H^0(\Lambda\oplus\Lambda^*)$.
Then
\beq{kurF}
\sigma_2\=\(p^*(t),-\tilde d^*\)\ \in\ H^0(\Lambda^*)\ \text{ cuts out }\F\text{ from }\tot_A(E_1^*).
\eeq
Therefore using $\Lambda$ to define the virtual cycle, the construction of \cite[Section 3.2]{OT1} gives the following (but see Remark \ref{but} below). The virtual cycle of $T^*[-2]\;\M$ is made by taking the intersection of 
$$
C_{F/\mathrm{tot}\_{\!A}(E_1^*)}\ \subset\ \Lambda^*|_F\ \text{ with the zero section }\ 0_{\Lambda^*|_F}\,,
$$
cosection localised by $\sigma_1|_F\;$:
\beq{virloc}
\pm\,0_{\Lambda^*|\_F}^{\,!,\,\loc}\big[C_{F/\mathrm{tot}\_{\!A}(E_1^*)}\big]\ \in\ A_{\;p-e}(M).
\eeq
Here we think of $\sigma_1|_F$ as a function on $\Lambda^*|_F$ (linear on the fibres) which vanishes identically on $C_{F/\mathrm{tot}\_{\!A}(E_1^*)}$ by \cite[Lemma 3.1]{OT1}. Thus the cosection localisation of \cite{KLi} applies, localising the intersection to the zeros of $\sigma_1|_F$ on the zero locus of $\sigma_2$, i.e. to the zero locus $M$ of $\sigma$ as claimed.

We need to describe the sign $\pm$ \eqref{virloc}, written in  \cite[Section 3.2]{OT1} as $(-1)^{|\Lambda|+\rk\Lambda}$. Recall we gave $K\oplus K^*$ the orientation $o_K$ for which $K$ is a positive maximal isotropic. Under the passage from \eqref{TT} to \eqref{t*t} this corresponds to giving
$p^*(B\oplus E_0\oplus E_0^*\oplus B^*)$ the orientation $o_{K\oplus p^*E_1^*}=o_K\otimes o_{p^*E_1^*}$ by \cite[Equation (65)]{OT1}. Writing this as $(-1)^{|\Lambda|}o_{\Lambda}$ defines $(-1)^{|\Lambda|}$.

Working locally we may split the exact sequence $0\to K\to p^*(B\oplus E_0)\to p^* E_1\to0$. Then suppressing some $p^*$s and setting $b=\rk B$, etc,
\begin{eqnarray}\nonumber
o_K\otimes o_{E_1^*} &=& (-1)^{e_1}o_K\otimes o_{E_1}\=(-1)^{e_1}o_{K\oplus E_1}\=(-1)^{e_1}o_{B\oplus E_0} \\
&=& (-1)^{e_1}o_{B}\otimes o_{E_0}=(-1)^{e_1+b}o_{B^*}\otimes o_{E_0}=(-1)^{e_1+b}o_{\Lambda}. \label{oo}
\end{eqnarray}
Thus our sign $\pm$ is $(-1)^{|\Lambda|+\rk\Lambda}=(-1)^{e_1+b+b+e_0}=(-1)^e$ as required.

Finally, by \eqref{sig} we see that $\sigma_1|_F=(p^*(s_0),\,\tilde s_1^*)\colon p^*(B\oplus E_0^*)|_F\to\cO_F$, which is the cosection $h^1(d\tilde s)$ as calculated in \eqref{cosex}. Thus \eqref{virloc} is precisely $(-1)^e\;[\F]^{\loc}$, as required.

\begin{rmk}\label{but}
More precisely, \eqref{virloc} is in fact $(-1)^{|\Lambda|}\sqrt e\,(\Lambda\oplus\Lambda^*,\sigma,\Lambda)$ \cite[Section 3.2]{OT1}\,---\,the $\sigma$-localised square root Euler class of $\Lambda\oplus\Lambda^*$. Its construction uses the family of isotropic graphs $\Gamma_{(\sigma_1,\,z^{-1}\sigma_2)}$ in $\Lambda\oplus\Lambda^*$ to interpolate between $\Gamma_{\sigma}$ (at $z=1$) and its degeneration $C_{F/\mathrm{tot}\_{\!A}(E_1^*)}$ (at $z=0$)\,---\,see the proof of \cite[Lemma 3.1]{OT1}. In contrast, the virtual cycle is defined in \cite[Section 4.2]{OT1} by replacing the graph by $C_{M/\mathrm{tot}\_{\!A}(E_1^*)}$. But this is also a degeneration of $\Gamma_{\sigma}$ through the family of isotropic graphs $\Gamma_{\!z^{-1}\sigma}$\,---\,see \cite[Equation (71)]{OT1}\,---\,so the result is the same.

We can combine these two families\footnote{The two families lie over $(w=1)\subset\C^2$ and $(z=1)\subset\C^2$ respectively. References for this section are \cite{BCM, KKP, Ma}.} over $\C$ into a single family over $\C^2$ by considering the (closure of the) isotropic graph
$$
\overline\Gamma\ :=\ \overline{\Gamma_{(w^{-1}\sigma_1,\,(wz)^{-1}\sigma_2)}}\ \subset\ 
\mathrm{tot}\_{\,\mathrm{tot}\_{\!A}(E_1^*)\times\C^2}
\big(\Lambda\oplus\Lambda^*\big).
$$
Here and below $w,z$ are the coordinates pulled back from $\C^2$ and we suppress some pullback maps on $\Lambda\oplus\Lambda^*$.
Denote the inclusion of the point $(w,z)$ into $\C^2$ by $i_{w,z}$. Then factoring $i_{1,0}$ and $i_{0,0}$ through the inclusion of $(z=0)$ gives a rational equivalence
$$
i_{1,0}^{\,!}\,\overline\Gamma\ \sim\ i_{0,0}^{\,!}\,\overline\Gamma\ \text{ inside }\ \overline\Gamma\times\_{\C^2}(z=0)\ \subset\ \mathrm{tot}\_{F\times(z=0)}\big(\Lambda\oplus\Lambda^*\big).
$$
Similarly factoring $i_{0,1}$ and $i_{0,0}$ through the inclusion of $(w=0)$ gives a rational equivalence
\beq{fam}
i_{0,1}^{\,!}\,\overline\Gamma\ \sim\ i_{0,0}^{\,!}\,\overline\Gamma\ \text{ inside }\ \overline\Gamma\times\_{\C^2}(w=0)\ \subset\ \mathrm{tot}\_{F\times(w=0)}\big(\Lambda\oplus\Lambda^*\big).
\eeq
But $i_{1,0}^{\,!}\,\overline\Gamma=C_{F/\mathrm{tot}\_{\!A}(E_1^*)}$ because $\overline\Gamma|_{(w=1)}$ is the (flat!) deformation of tot$\_{\!A}(E_1^*)$ to the normal cone $C_{F/\mathrm{tot}\_{\!A}(E_1^*)}$. Similarly 
$i_{0,1}^{\,!}\,\overline\Gamma\=C_{M/\mathrm{tot}\_{\!A}(E_1^*)}$ because $\overline\Gamma|_{(z=1)}$ is the deformation of tot$\_{\!A}(E_1^*)$ to the normal cone $C_{M/\mathrm{tot}\_{\!A}(E_1^*)}$.
The upshot is an \emph{isotropic} rational equivalence
\beq{gam}
C_{F/\mathrm{tot}\_{\!A}(E_1^*)}\ \sim\ C_{M/\mathrm{tot}\_{\!A}(E_1^*)}
\ \text{ inside }\ \mathrm{tot}\_{F\times(wz=0)}\big(\Lambda\oplus\Lambda^*\big).
\eeq
This will prove useful later because it replaces the ambient space $\mathrm{tot}\_{\!A}(E_1^*)$ (which only exists locally in general) with data over $F$ (which will globalise). We note that the embedding of the base $F$ of $C_F$ into $\Lambda\oplus\Lambda^*$ factors through the first factor as $\sigma_1=(\pi^*s_0,\,\tilde s_1^*)=h^1(d\;\tilde s)$. 
\end{rmk}

\section{Proof of main result}\label{last}
\emph{Throughout this section we always restrict to the underlying scheme or stack $Y=\pi_0({\bf Y})$ of whichever derived space $\bf Y$ we are working on.} We usually omit the restriction map for brevity.

Since all of our stacks have the resolution property we may always resolve (complexes of) sheaves using complexes of very negative locally free sheaves. So to prove \eqref{thm} in the general case we begin by finding global resolutions of the objects $\TT_{\P},\,\TT_{\F},\,\TT_{\M},\,E\udot$ and the maps between them. We will choose them to be reminiscent of the (derivative of) the local model \eqref{localmodel}.\footnote{The $A,B,t',\cA,E_0,E_1,s_0,s_1$ of this Section play the roles of the $T_A,B,dt,T_{\tot\_{\!A}(E_1^*)}$, $E_0,E_1,s_0,s_1$ (all restricted to $P,\,F$ or $M$) of Sections \ref{cosec} and \ref{idea}.}\smallskip

First pick a 2-term locally free resolution $E_0\rt{d}E_1$ of $\Ed$. Now we may pick a 2-term locally free resolution $t'\colon A\to B$ of $\TT_{\P}$ with $B$ sufficiently negative that the functor $\Ext^{>0}(B,\ \cdot\ )$ is zero on $E_0^*,\,E_1^*$ and $\iota_*(E_0^*|_M\)$.

\subsection*{Representative for $\TT_{\F}$}
The exact triangle $\TT_{\F}\to\pi^*\TT_{\P}\to\TT_{\F/\P}[1]=\pi^*\Ed^\vee$ defines an element of the uppermost group in the diagram
$$
\xymatrix@R=15pt{
&& \Hom(\pi^*\TT_{\P},\pi^*\Ed^\vee) \ar[d] \\
0 \ar[r]& \Ext^1(\pi^*A,\pi^*E_1^*) \ar[d]\ar[r]^-\sim& \Ext^1(\pi^*\TT_{\P},\pi^*E_1^*)  \ar[r]\ar[d]& 0 \\
0 \ar[r]& \Ext^1(\pi^*A,\pi^*E_0^*) \ar[r]^-\sim& \Ext^1(\pi^*\TT_{\P},\pi^*E_0^*)  \ar[r]& 0.\!\!}
$$
Here the horizontal arrows come from the triangle $\TT_{\P}\to A\rt{t'}B$ and the $\Ext^{>0}(B,\ \cdot\ )$ vanishing. Thus we get an extension in $\Ext^1(\pi^*A,\pi^*E_1^*)$ which maps to zero in $\Ext^1(\pi^*A,\pi^*E_0^*)$, inducing a commutative diagram
$$
\xymatrix@R=15pt@C=30pt{
\pi^*E_1^* \ar[r]\ar[d]_{-\pi^*d^*}& \cA \ar[r]\ar[d]& \pi^*A \ar@{=}[d] \\
\pi^*E_0^* \ar[r]^-{(0,1)}& \pi^*(A\oplus E_0^*) \ar[r]^-{(1,0)}& \pi^*A.\!}
$$
Composing with $\pi^*t'\colon\pi^*A\to\pi^*B$ gives the representatives
\beq{AA}
\xymatrix@R=0pt@C=10pt{
\pi^*E_1^* \ar[rr]\ar[dd]_{-\pi^*d^*}&& \cA \ar[rr]\ar[dd]&& \pi^*A \ar[dd]^{\pi^*t'} \\
&&&&& \text{for}\hspace{-3mm} & \pi^*\Ed^\vee[-1] \ar[r]+<-6pt,0pt>& \ \ \TT_{\F} \ar[r]+<-15pt,0pt>& \ \pi^*\TT_{\P} \\
\pi^*E_0^* \ar[rr]^-{(0,1)}&& \pi^*(B\oplus E_0^*) \ar[rr]^-{(1,0)}&& \pi^*B}
\eeq
because the connecting homomorphism of the horizontal triangle of vertical 2-term complexes represents $\pi^*\TT_{\P}\to\TT_{\F/\P}[1]=\pi^*\Ed^\vee$ by construction.

The zero section $\P\subset\F$ defines a splitting $\TT_{\P}\to\TT_{\F}|_P$ of the triangle $\TT_{\F}\to\pi^*\TT_{\P}\to\TT_{\F/\P}[1]$ on $P$. We note for later that following through the above construction then shows that on restriction to $P\subset F$, \eqref{AA} splits as
\beq{AB}
\xymatrix@R=18pt{
E_1^* \ar[r]^-{(1,0)}\ar[d]_{-d^*}& E_1^*\oplus A \ar[r]^-{(0,1)}\ar[d]<.5ex>_{-d^*\oplus\!}^{t'}& A \ar[d]^{t'} \\
E_0^* \ar[r]^-{(1,0)}& E_0^*\oplus B \ar[r]^-{(0,1)}& B.\!\!}
\eeq

\subsection*{Representative for $\TT_{\M}$} Consider the diagram of exact triangles
\beq{ABTP}
\xymatrix@R=18pt{
B|_M[-1] \ar[r]\ar@{..>}[d]^{s_1}& \TT_{\P}|_M \ar[r]\ar[d]^{ds}& A|_M \ar@{..>}[d]^{ds_0} \\
E_1|_M[-1] \ar[r]& \Ed|_M \ar[r]& E_0|_M.\!\!}
\eeq
Since $\Ext^1(B|_M,E_0|_M)=0$ there exists a map $s_1$; taking cones then gives the map marked $ds_0$. (Due to the choices involved $ds_0$ is not entirely determined by $s_0$\,---\,the composition of $s$ with $\Ed\to E_0$\,---\,so the notation is only suggestive.) Hence we get a representative of $\TT_{\M}$\,---\,the cocone of $\TT_{\P}|_M\to\Ed^\vee|_M$\,---\,as the total complex of
$$
\xymatrix@R=18pt{
A|_M \ar[r]^{t'}\ar[d]_{ds_0} & B|_M \ar[d]^{s_1} \\
E_0|_M \ar[r]^d& E_1|_M.\!}
$$
Since $\M$ is quasi-smooth $(s_1,-d)\colon(B\oplus E_0)|_M\to E_1|_M$ is onto. Letting $K$ denote its kernel we get three quasi-isomorphic representatives of $\TT_{\M}$,
\begin{eqnarray}
\TT_{\M} &\cong& \big\{A|_M\rt{(t',ds_0)}(B\oplus E_0)|_M\rt{(s_1,-d)}E_1|_M\big\} \label{TM1} \\
&\cong& \big\{A|_M\rt{(t',ds_0)}K\big\} \label{TM2} \\
&\cong& \big\{(A\oplus E_1^*)|_M\rt{(t',ds_0)\oplus\id}K\oplus E_1^*|_M\big\}. \label{TM3}
\end{eqnarray}
Since $\TT^\vee_{\M}|\_M\to\LL_M$ is a perfect obstruction theory, the Behrend-Fantechi construction \cite{BF} applied to the third complex \eqref{TM3} defines a cone
\beq{CM}
C_M\ \subset\ K\oplus E_1^*|_M \quad\text{such that}\quad [\M]^{\vir}\=0^{\;!}_{K\oplus E_1^*|_M}[C_M].
\eeq

For later we note that by \eqref{ABTP} the map $\TT_{\M}|_M\to\TT_{\P}|_M$ becomes the projection of the complex \eqref{TM1} to its first two terms $t'\colon A|_M\to B|_M$. Thus it is also represented by the chain map from \eqref{TM2} to $A|_M\to B|_M$ given by the identity on $A|_M$ and the composition $K\into(B\oplus E_0)|_M\to B|_M$ on $K$.

We further compose this with the (restriction to $M$ of the) map $\TT_{\P}\to\TT_{\F}|_P$ induced by the zero section $\P\subset\F$, to get a description of the map
\beq{MF}
\TT_{\M}|_M\To\TT_{\F}|_M
\eeq
induced by $\M\subset\F$. By \eqref{AB} $\TT_{\P}\to\TT_{\F}|_P$ is the inclusion of $t'\colon A\to B$ as the second summand of the central vertical complex in \eqref{AB}. Thus \eqref{MF} maps \eqref{TM2} in the obvious way to the second summand of
\beq{MF2}
\xymatrix@R=-4pt{
& E_1^*|\_M & E_0^*|\_M \\
\TT_{\F} \qquad\cong& \quad\oplus\quad \ar[r]^{-d^*}_{\oplus t'}& \quad\oplus\quad \\
& A|\_M & B|\_M.\!\!}
\eeq

\subsection*{Representative for $\TT_{T^*[-2]\M}$}
On restriction to $M$ we have the splitting\footnote{The derivative of the zero section $\M\to T^*[-2]\M$ splits the pullback to the zero section of the exact triangle $\TT_{\M}^{\vee}[-2]\to\TT_{T^*[-2]\M}\to\TT_{\M}$.} $\TT_{T^*[-2]\M}\cong\TT_{\M}\oplus\TT^{\vee}_{\M}[-2]$, represented by the complex $\eqref{TM1}\;\oplus\;\eqref{TM1}^\vee[-2]$,
\beq{t*}
\quad (A\oplus E_1^*)\big|_M\xrightarrow[(-d^*\!,\,s_1^*)]{(t'\!,\,ds_0)\oplus}(B\oplus E_0)\big|_M\oplus(E_0^*\oplus B^*)\big|_M\To(A^*\oplus E_1)\big|_M.
\eeq
Here the second arrow is the dual of the first. In fact \eqref{t*} is also equal to $\eqref{TM3}\;\oplus\;\eqref{TM3}^\vee[-2]$ because $K\into(B\oplus E_0)|_M$ induces $E_1^*\cong K^\perp\into(E_0^*\oplus B^*)|_M$ and hence an isomorphism between $(K\oplus K^\perp)\oplus(K\oplus K^\perp)^*$ and $(B\oplus E_0\oplus E_0^*\oplus B^*)|_M$.

Thus the first arrow of \eqref{t*} factors through the maximal isotropic subbundle $K\oplus E_1^*|_M$. So when we consider the stupid truncation of \eqref{t*}
\beq{trunc}
\TT_{\tau\M}\ :=\ \big\{(A\oplus E_1^*)\big|_M\To(B\oplus E_0)\big|_M\oplus(E_0^*\oplus B^*)\big|_M\big\}
\eeq
as defining a perfect obstruction theory $\TT_{\tau\M}^\vee\big|_M\to\LL_M$ for $M$, the induced Behrend-Fantechi cone is
\beq{CM*}
C_M\,\stackrel{\eqref{CM}}{\subseteq}\,K\oplus E_1^*|_M\ \subset\ (B\oplus E_0\oplus E_0^*\oplus B^*)|_M.
\eeq
The virtual cycle of $T^*[-2]\;\M$ is defined in \cite[Section 3.3]{OT1}\footnote{This definition requires a choice of orientation on $(B\oplus E_0\oplus E_0^*\oplus B^*)|_M$. We use the choice $o_{K\oplus E_1^*|_M}$ which makes $K\oplus E_1^*|_M$ a \emph{positive} maximal isotropic.} via this stupid truncation as
\beq{a}
\sqrt{0^{\,!}}_{(B\oplus E_0\oplus E_0^*\oplus B^*)|_M}[C_M]\=0^{\,!}_{K\oplus E_1^*|_M}[C_M]\ \stackrel{\eqref{CM}}{=}\ [\M]^{\vir},
\eeq
where the first equality is \cite[Lemma 3.5]{OT1} applied to the isotropic embedding \eqref{CM*}.

\subsection*{Relating $\M$ and $\F$} The stupid (dual) perfect obstruction theory \eqref{trunc} sits inside the exact triangle
$$
\xymatrix@R=15pt@C=18pt{
(E_0\oplus B^*)\big|_M[-1]\hspace{-5mm} \ar[d]&&& \hspace{18mm}(E_0\oplus B^*)\big|_M \ar[d]<4ex> \\
\TT_{\tau\M} \ar[d] &=& (A\oplus E_1^*)\big|_M \ar@{=}[d]<-1ex>\ar[r]& (B\oplus E^*_0\oplus E_0\oplus B^*)\big|_M \ar[d]<-7.5ex> \\
\TT_{\F}\big|_M && (A\oplus E_1^*)\big|_M \ar[r]& (B\oplus E_0^*)\big|_M.\hspace{18mm}}\vspace{-22mm}
$$
\beq{ton}\vspace{13mm}\eeq
On the bottom row we have used the splitting \eqref{AB}, which also shows the obvious vertical arrows are chain maps.

We next show this fits into the commutative diagram of perfect obstruction theories \eqref{tri} below. By (\ref{MF}, \ref{MF2}) $\TT_{\M}\to\TT_{\F}|_M$ factors through the above map $\TT_{\tau\M}\to\TT_{\F}|_M$ by the natural inclusion of \eqref{TM2} into \eqref{trunc}. This defines the top row of the commutative diagram
$$
\xymatrix@R=16pt{
\TT^\vee_{\F}\big|_M \ar[r]\ar[d]& (\TT_{\tau\M})^\vee \ar[r]\ar[d]&
\TT^\vee_{\M} \ar[d] \\
\LL_F\big|_M \ar[r]& \LL_M \ar@{=}[r]& \LL_M,\!\!}
$$
with the vertical maps induced by $F\into\F$ and $M\into\M$ respectively.
Combining the left hand square with the dual of \eqref{ton} gives the exact triangle of perfect obstruction theories
\beq{tri}
\xymatrix@R=16pt{
(E_0^*\oplus B)\big|_M \ar[r]^-{(ds_0)^*}_-{\oplus s_1}\ar[d]& \TT^\vee_{\F}\big|_M \ar[r]\ar[d]& (\TT_{\tau\M})^\vee \ar[r]\ar[d]&
(E_0^*\oplus B)\big|_M[1] \ar[d] \\
\LL_{M/F}[-1] \ar[r]& \LL_F\big|_M \ar[r]& \LL_M \ar[r]& \LL_{M/F}.\!\!}
\eeq

\subsection*{Isotropic rational equivalence of cones}
Applying the Behrend-Fantechi construction to the perfect obstruction theory \eqref{AA} for $F$ gives a cone
\beq{cf}
C_F\ \subset\ \pi^*(B\oplus E_0^*)\ =:\ \Lambda^*\ \subset\ \Lambda\oplus\Lambda^*.
\eeq
Then by \cite{KKP} the exact triangle \eqref{tri}\,---\,and its realisation \eqref{ton} as a short exact sequence of 2-term chain complexes of locally free sheaves\,---\,induces a rational equivalence between
$$
C_M\ \subset\ (\Lambda\oplus\Lambda^*)\big|_M\ \text{ and }\ 
C_{M/C_F}\ \subset\ (\Lambda\oplus\Lambda^*)\big|_M\,.
$$
In order to prove this rational equivalence takes places inside an isotropic substack (in fact cone) of $(\Lambda\oplus\Lambda^*)\big|_M$ we review parts of the \cite{KKP} construction. For more details see also \cite{BCM, Ma}.

Fix any Deligne-Mumford stacks $M\subset F$ with perfect obstruction theories $\TT_M^\vee\to\LL_M,\ \TT_F^\vee\to\LL_F$ fitting into a diagram of exact triangles
$$
\xymatrix@R=15pt{
\Lambda[-1] \ar[r]& \TT_M \ar[r]& \TT_F|_M \ar[r]& \Lambda \\
\LL_{M/F}^\vee \ar[r]\ar[u]& \LL_M^\vee \ar[r]\ar[u]& \LL_F^\vee|_M \ar[r]\ar[u]& \LL_{M/F}^\vee[1],\!\! \ar[u]}
$$
with the (first three terms of the) top row represented by a short exact sequence of (vertical) 2-term complexes of vector bundles
\beq{kkpee}
\xymatrix@R=15pt{
& \cA \ar@{=}[r]\ar[d]& \cA\ar[d] \\
\Lambda\ \ar@{^(->}[r]<-.2ex>& E \ar@{->>}[r]^\pi& \Lambda^*.\!\!}
\eeq
For us these will be provided by \eqref{ton} and \eqref{tri} with $\cA=(A\oplus E_1^*)|_M,\ \Lambda=(B^*\oplus E_0)|_M$ and $E=\Lambda\oplus\Lambda^*$.

Using this data \cite[Proposition 1]{KKP} defines a canonical abelian cone stack (a certain normal sheaf) inside a bundle stack
$$
[N/\cA]\ \subset\ [E_t/\cA]\,\To\,M\times\C.
$$
Here $E_t$ is the bundle over $M\times\C$ given by degenerating the extension $E$ to its splitting $\Lambda\oplus\Lambda^*$; over $t\in\C$ it is the kernel of $(t\id,\,\pi)\colon\Lambda^*\oplus E\onto\Lambda^*$. (In our situation $E=\Lambda\oplus\Lambda^*$ is already split so $E_t$ is just the pullback of $E$ from $M$ to $M\times\C$.) We have also pulled $\cA$ back to $M\times\C$.

Then \cite[Equation (8)]{KKP} defines a canonical normal cone stack $[C/\cA]\subset[N/\cA]$ containing a rational equivalence between any $t\ne0$ fibre $[C_M/\cA]$\,---\,the intrinsic normal cone of $M$\,---\,and the central fibre $[C_{M/C_F}/\cA]$ over $t=0$. Pulling back by $E_t\to[E_t/\cA]$ gives, in our situation, a canonical cone $C\subset\mathrm{tot}_{M\times\C}(\Lambda\oplus\Lambda^*)$ satisfying
\begin{itemize}
\item the fibre of $C$ over points $i_t\colon\{t\}\into\C$ with $t\ne0$ is $C_M$,
\item the fibre of $C$ over $t=0$ \emph{contains} $C_{M/C_F}$, and
\item $i_0^{\;!}[C]=[C_{M/C_F}]$.
\end{itemize}
This gives our rational equivalence between $i_1^{\;!}[C]=[C_M]$ and $[C_{M/C_F}]$ inside $C\subset(\Lambda\oplus\Lambda^*)|_M$.

About any point of $M$ our (exact triangle of) perfect obstruction theories (\ref{ton}, \ref{tri}) is isomorphic to one arising from the local model \eqref{localmodel}.\footnote{After shrinking $M$ and $P$ we can find a local Kuranishi structure $(A,B,t)$ for $P$ which induces the perfect obstruction theory $dt|_P=t'\colon A|_P\to B|_P$; see \cite[Theorem 3.3]{OT2} for example. (Here $A$ may be taken to be an open set in a vector space, excusing our abuse of notation in identifying it with its tangent spaces $A$. We are also using $B$ to denote both the bundle on $P$ and a choice of a local extension of it to $A$.) Then pick local lifts of $d\colon E_0\to E_1$ and $s$ to $A$ and proceed as in \eqref{localmodel} and Section \ref{idea} in the ambient space $\cA=\tot_A(E_1^*)\rt pA$. This gives Kuranishi charts (local over $M$ but global in the $p$-fibre directions) $(\cA,\Lambda\oplus\Lambda^*,\sigma)$ for $\tau\;\bf M$ \eqref{Darboux} and $(\cA,\Lambda^*,\sigma_2)$ for $\bf F$ \eqref{kurF}, compatible under the projection $\Lambda\oplus\Lambda^*\to\Lambda^*$ which maps $\sigma=(\sigma_1,\sigma_2)$ \eqref{sig} to $\sigma_2$ \eqref{kurF}. Taking their derivatives along $M$ to pass back to perfect obstruction theories recovers precisely (\ref{ton}, \ref{tri}) and so the exact triangle \eqref{kkpee} required to apply \cite{KKP}. In this local model \cite{KKP}'s cone $C$ is, by construction, the cone $\overline\Gamma\times\_{\C^2}(w=0)$ \eqref{fam} of Remark \ref{but}.}

In this local model $C_M$\,---\,the limit of the isotropic graphs \eqref{Darboux}\,---\,and $C$\,---\,the family over $(w=0)$ in Remark \ref{but}\,---\,are \emph{isotropic} in $\mathrm{tot}_{M\times\C}(\Lambda\oplus\Lambda^*)$. Since $C$ is canonical it is isomorphic to its local model and is therefore also isotropic.

Therefore we can apply the deformation invariance \cite[Equation (78)]{OT1} of $\surd\;0^{\;!}$ to \eqref{a} to give
\beq{nearly}
[\M]^{\vir}\ =\ (-1)^{|\Lambda|+\rk\Lambda}\sqrt{0^{\,!}_{\Lambda\oplus\Lambda^*}}\,[C_{M/C_F}].
\eeq

Finally we want to replace $C_{M/C_F}$ by $C_F$ by deforming the embedding $(0,c)\colon C_F\into\Lambda\oplus\Lambda^*$ of \eqref{cf}, where $c$ is the embedding $C_F\into\Lambda^*$.

Writing $\TT_{\F}=\{\cA\to\Lambda^*\}$ as in \eqref{AA} we can consider the composition
$$
\Lambda^*\To h^1(\TT_{\F})\rt{h^1(d\;\tilde s)}\cO_F
$$
as a section of $\Lambda$ which we also denote by $h^1(d\;\tilde s)$. We use this to perturb $(0,c)$ \eqref{cf}, taking the closure of the graph
\beq{closure}
\overline{\Gamma_{(w^{-1}h^1(d\;\tilde s),\,c)}}\ \subset\ \mathrm{tot}\_{F\times\PP^1}(\Lambda\oplus\Lambda^*),
\eeq
where $w$ is the coordinate pulled back from $\PP^1$. Because $h^1(d\;\tilde s)$ cuts out $M\subset C_F$ this gives the standard (flat!) deformation of $C_F$ to the normal cone of $M\subset C_F$. It therefore gives a rational equivalence between $C_{M/C_F}\subset(\Lambda\oplus\Lambda^*)|_M$ and the fibre \eqref{cf} over $w=\infty$. 

Since $h^1(d\;\tilde s))$ is a cosection the composition
$$
C_F\Into\Lambda^*\rt{h^1(d\;\tilde s)}\cO_F
$$
vanishes, which means the rational equivalence \eqref{closure} is \emph{isotropic}.
So by the deformation invariance \cite[Equation (78)]{OT1} again, \eqref{nearly} has become
$$
[\M]^{\vir}\ =\ (-1)^{|\Lambda|+\rk\Lambda}\sqrt{0^{\,!}_{\Lambda\oplus\Lambda^*}}\;[C_F],
$$
where $C_F$ is embedded in $\Lambda\oplus\Lambda^*$ via $(h^1(d\;\tilde s),c)$.
%
%
Finally, this class is defined in \cite[Section 3.2]{OT1} to be the intersection of $C_F$ with the 0-section of $\Lambda^*$, cosection localised by the tautological cosection $\tau\_{\Lambda}\big|_F$ of the pullback of $\Lambda$ to $F$. But since $F$ is embedded in this pullback by the graph of $h^1(d\;\tilde s)$, this cosection is just $h^1(d\;\tilde s)$, yielding
$$
[\M]^{\vir}\ =\ (-1)^{|\Lambda|+\rk\Lambda}\,0^{\,!,\,\mathrm{loc}}_{\Lambda^*,\,h^1(d\;\tilde s)}[C_F]\ =\ (-1)^e[\F]^{\mathrm{loc}}.
$$
The verification of the sign $(-1)^{|\Lambda|+\rk\Lambda}=(-1)^e$ was done in \eqref{oo}.

\bibliographystyle{halphanum}
\bibliography{references}

\bigskip \noindent {\tt{j.oh@imperial.ac.uk \\ richard.thomas@imperial.ac.uk}} \medskip

\noindent Department of Mathematics \\
\noindent Imperial College London\\
\noindent London SW7 2AZ \\
\noindent United Kingdom

\end{document}